\documentstyle{amsppt}
\magnification=1100
\input xy
\xyoption{all}

\topmatter
\title On the 0-dimensional irreducible components of the singular locus of $\Cal A_g$  \endtitle
\rightheadtext{On the singular locus of $\Cal A_g$} \subjclass
14K10, 14K22 \endsubjclass
\author V. Gonz\'alez -Aguilera, J.M. Mu\~ noz-Porras  and Alexis G. Zamora \endauthor
\address Departamento de Matem\'aticas, Universidad Santa Mar\'{\i}a,\newline
Avenida Espa\~ na 1680
\newline Valpara\'{\i}so, Chile.\endaddress
\email victor.gonzalez\@ usm.cl \endemail
\address Departamento de Matem\'aticas, Universidad de Salamanca\newline
Plaza de la Merced 1-4,
\newline Salamanca 37008, Espa\~ na \endaddress
\email jmp\@ usal.es \endemail
\address CIMAT,\newline Callej\'on de Xalisco s/n
\newline Valenciana, Guanajuato, Gto.,\newline C.P. 36240, \newline M\'exico \endaddress
\email alexis\@ cimat.mx \endemail \abstract In this short note we
give a characterization of extremal principally polarized abelian
varieties determining an isolated point in $Sing \Cal A_g$. The
case $g=5$ is treated with detail.
\endabstract
\thanks The first author was partially supported by Fondecyt Grant
N1010432 and DGIP UTFSM Grant 120321. The second author was
partially supported by MCT Grant BFM 2003-00078 and JCyL Grant
SA071/04. The third author was partially supported by Conacyt
Grant 41459-F.
\endthanks
\endtopmatter

\document

\head 1. Introduction \endhead

Let $\Cal A_g$ be the moduli space of principally polarized
abelian varieties (p.p.a.v.). In [3] it was proved the existence
of irreducible varieties $\Cal A_g(p,\rho_n)\subset \Cal A_g$
representing triples $(X,\sigma, \alpha)$ consisting of a p.p.a.v.
$X$, a non trivial automorphism $\sigma \in Aut(X)$ of order $p$
(a prime number) and the conjugation class $\alpha$ of the
representation of $\sigma$ in $X_n$ for some fixed $n\ge 3$. The
local information about these varieties is encoded in the so
called "analytic representation" of $\sigma$, that is, the
differential

$$ d\sigma: T_0X @>>> T_0X.$$

Indeed, if $\{ \xi^{i_1},...,\xi^{i_g}\}$ denotes the set of
eigenvalues of $d\sigma$ then

$$\dim \Cal A_g(p,\rho_n)= \frac{n_0(n_0+1)}{2}
+\sum_{i=1}^{p-1/2}n_in_{p-i},$$ where $n_i$ denotes the dimension
of the eigenspace associated to $\xi^i$. Moreover, there exists a
number $r$ such that $n_i+n_{p-i}=r$ for $i\ne 0$ and

$$g=n_0 +\frac{p-1}{2}r.$$ The existence of such number is a consequence of
the fact that the rational representation of $\sigma$ is the
direct sum $d\sigma\oplus \bar{d\sigma}$ ([3], Prop.2.1).

Two fundamental problems arise from these results. In the first
place there is a problem of existence: given numerical data
satisfying the above restriction, it there exists an p.p.a.v. $X$
admitting an automorphism compatible with this data? In the second
place, in order to describe the irreducible components of $Sing
\Cal A_g$ we need to determine what inclusions of the form $\Cal
A_g(p,\rho_n) \subset \Cal A_g(q,\rho'_n)$ occurs. As we shall
show both problems are closely related.

The first was solved positively in [3] in the case $n_0=0$, the
idea is to reduce the existence problem to the existence of
p.p.a.v. of dimension $g$ admitting an automorphism of prime order
$=2g+1$ (we shall call such a p.p.a.v. an extremal one), note that
in this case $r=1$, $n_0=0$ and $n_in_{p-i}=0$. The existence of
these varieties is ensured by a classical construction ([5],
section 22, [L], Chapter 1) which in fact is valid in a more
general context: the existence of CM-type Abelian Varieties. Then
the existence problem is solved by considering products of these
"extremal varieties" ([3], Lemma 2.7).

In this paper we study and solve the problem of inclusion of
varieties $\Cal A_g(p,\rho_n) \subset \Cal A_g(q,\rho'_n)$ in the
extremal case, i.e. under the assumption $p=2g+1$. It should be
noticed that if the number $2g+1$ is a prime, then it is the
maximal possible order for an automorphism of prime order acting
on p.p.a.v. of dimension $g$. It is clear that $\dim \Cal
A_g(2g+1, \rho)=0$.

In order to fix the notation we recall the construction of
extremal p.p.a.v. The starting point is $\Bbb F_p^*$ which we
identify with $\{\xi,\xi^2,...,\xi^{p-1}\}$, next, fix a partition
$\{C,C'\}$ of $\Bbb F_p^*$ with $\vert C\vert=g$. Let
$C=\{\xi^{i_1},...,\xi^{i_g}\}$. Define the isomorphism:

$$\Phi_C: \Bbb Q(\xi)\otimes_{\Bbb Q} \Bbb R @>>> \Bbb C^g,$$ given
by $\Phi_C(\xi)=(\xi^{i_1},...,\xi^{i_g})$. Obviously $\xi$ acts
on $\Bbb C^g$ as a linear automorphism of order $p$ under this
identification. Let $I$ be a fractional ideal of $\Bbb Q(\xi)$.
$L:=\Phi_C(I)\subset \Bbb C^g$ is a lattice invariant under the
action of $\xi$. In this way $X=\Bbb C^g/L$ is a complex torus
admitting an action of $\Bbb Z_p$ with analytic representation
determined by the set of eigenvalues $C$. A Riemann form can be
defined on $(\Bbb C^g, L)$ in such a way that $X$ becomes an
extremal p.p.a.v.

Multiplication defines a natural action of $\Bbb F_p^*$ in the set
$\Cal C$ of all the subsets of $\Bbb F_p^*$ of cardinality $g$.
Let $H_C$ be the isotropy group of $C$ under this action. The main
result of this note is the following:

\proclaim{Theorem 1}. Let $X$ be an extremal p.p.a.v.. Denote by
$\sigma$ the automorphism of $X$ of order $2g+1$. Assume that the
 analytic representation $d\sigma$ has as set of eigenvalues the
 set $C$. Then $X$ is an isolated point in $Sing \Cal A_g$ if and only
 if $H_C$ is trivial.\endproclaim

 Section 2 of the paper is devoted to prove this theorem. After
 this we include in section 3 a classification of the extremal
 p.p.a.v. of dimension $5$ that determine isolated points in $Sing
 A_5$. The proof of Theorem 1 is based in our previous paper [3]
 as well as in the results in [2].

 \head 2. Proof of Theorem 1 \endhead

\demo{Proof of Theorem 1}

Following the notation introduced in the previous section, let $X$
be an extremal p.p.a.v. constructed from the set $C$. $\Bbb F_p^*$
is the Galois group of $\Bbb Q(\xi)/\Bbb Q$, an element $k$ of
$\Bbb F_p^*$ acts in $\Bbb Q(\xi)$ by means of $k.\xi=\xi^k$. The
relevance of the group of isotropy $H_C\subset \Bbb F_p^*$ defined
in the introduction lies in the following observation: let $\theta
\in H_C$, $\theta \ne 1 $. $\theta$ can be extended to a $\Bbb
R$-linear automorphism:

 $$\theta:\Bbb R(\xi)@>>>\Bbb R(\xi),$$ as $\theta C=C$ the map
 $\Phi_C(\theta)$ is a $\Bbb C$-linear isomorphism. Furthermore, $\theta$
 determines an element of $Aut(X)$,
 that we are going to denote for simplicity just by $\theta$. In
 this way $H_C\ne 1$ implies that the total group of automorphism of
 $X$ contains two different cyclic subgroups: $<\xi>$ and
 $<\theta>$. Call $G$ to the total automorphism group, $G:=
 Aut(X)$, and $G_+:= G/\{\pm 1 \}$.

 We have  the following exact sequence of groups:

 $$1@>>> Z_G(<\xi>) @>>> N_G(<\xi>)@>\phi >> H_C @>>>1,\tag 1.1 $$ with
 $N_G(<\xi>)$ (respectively $Z_G(<\xi>)$ denoting as usual the
 normalizer (resp. the centralizer) of $<\xi>$ in $G$. The only
 map that needs definition is $\phi$, if $g\in N_G(<\xi>)$,
 $g.\xi.g^{-1}=\xi^k$, then we define $\phi(g)=k\in \Bbb F_p^*$. $\phi$ is surjective,
 as $\xi$ and $\xi^k$ are conjugate if and only if they have the
 same set of eigenvalues. Note
 that, by construction, $\theta \xi \theta^{-1}= \xi^{\theta}$.

 The sequence 1.1 induces the following:

 $$1@>>> <\xi> @>>> N_+@>\phi_+>> H_C @>>>1,\tag 1.2 $$
 where $N_+=N_G(<\xi>)/{\pm 1}$. In fact, $\phi_+(g)=1$ means
 $g\xi=\xi g$. Now, $g$ as an automorphism
 of $\Bbb Q(\xi)$ must acts as multiplication for some element: $g(x)=ux$,
 $\forall x\in \Bbb Q(\xi)$, and, as $g$ commutes with $\xi$,
 $uI=I$. It follows that $u\in \Bbb Q(\xi)^*$ and $g\in <\xi>$.

 Now, we can conclude the proof of the theorem: if $H_C=1$ then by
 the previous exact sequence $N_+=<\xi>$ and by Theorem 2 in [2]
 $G_+=<\xi>$.

 Conversely, assume $H_C\ne 1$ and let $\theta \in
 H_C$ be an element of order $q$, a divisor of $(p-1)$. The action of $\theta$ on $\Bbb Q(\xi)$
 gives rise to a decomposition:

 $$\Bbb Q(\xi)\simeq \Bbb Q(\xi)^{H_C}\oplus M,$$
 where $\Bbb Q(\xi)^{H_C}$ is the module of $H_C-$ invariants and $M$ its complement.
 Moreover, $\Bbb Q(\xi)^{H_C}\ne 0$,
 since $\theta (\xi^{k_1}+...\xi^{k_g})=\xi^{k_1}+...+\xi^{k_g}$. After tensoring
with $\Bbb R$ we obtain:

 $$V\simeq V_0\oplus ...\oplus V_m,$$

 with $V_0=\Bbb Q(\xi)^{H_C}\otimes \Bbb R\ne 0$. We
 conclude that the space of $\theta$-invariants $V_0$ is of
 strictly positive dimension. Thus, $X\in \Cal A_g(q,\theta_n)$ with

 $$\dim \Cal A_g(q,\theta_n)\ge 1.$$ \qed
 \enddemo

The previous characterization gives, moreover, a nice
classification of extremal p.p.a.v. in terms of the simplicity of
these varieties. This is provided by a theorem of Shimura and
Taniyama ([6], [4]). Explicitly:

\proclaim{Corollary 1.1} Let $X$ be an extremal p.p.a.v.. $X$ is
an isolated point of $Sing \Cal A_g$ if and only if $X$ is a
simple Abelian variety.\endproclaim

\demo{Proof} Just combine Theorem 1 with Theorem 3.5, Chapter 1 in
[4].\qed
\enddemo

Another useful remark is that there exists a simple sufficient
condition that guarantees the existence of extremal p.p.a.v. with
$H_C=1$. In fact, if $C=\{k_1,...,k_g\}$ and $\theta C=C$ then

$$k_1+...+k_g \cong \theta (k_1+...+k_g) \mod (p),$$
and, therefore, $k_1+...+k_g$ is divisible by $p$. In particular,
for each $g=(p-1)/2$, extremal Abelian varieties constructed from
$C=\{1,...,g\}$ are always isolated points in $Sing \Cal A_g$. The
argument is taken from [6], 8.4, page 72.

 \head 3. Extremal Abelian Varieties of dimension $5$ \endhead

 In this sections we work out with some detail the case $g=5$,
 $p=11$. First of all, the number of different (non-equivalent)
 partitions $C$ of $\Bbb F_{11}^*$ is $4$, they can be represented by:

 $$C_1=\{\xi,\xi^2,\xi^3,\xi^4,\xi^5\}, \quad
 C_2=\{\xi,\xi^2,\xi^3,\xi^4,\xi^6\},$$

 $$C_3=\{\xi,\xi^2,\xi^3,\xi^5,\xi^7\}, \quad
 C_4=\{\xi,\xi^3,\xi^4,\xi^5,\xi^9\}.$$

 Denote by $X_i$ the extremal p.p.a.v. corresponding to each of
 these partitions. These variety are uniquely determined by $C_i$,
 since $\Bbb Z [\xi]$ is a principal ideal domain if $\xi$ is a root of $11$. We have:

 \proclaim{Theorem 2} The only extremal p.p.a.v. in $\Cal A_5$ which
 are isolated in $Sing \Cal A_5$ are $X_1$, $X_2$ and $X_3$. Moreover:

 i) $X_1$ is isomorphic to the Jacobian variety of the
 hyperelliptic curve $y^2=\prod_{j=1}^{11}(x-\xi^j)$.

 ii) $X_2$ is isomorphic to the Jacobian variety of the
 Lefschetz curve $y^{11}=x^2(x-1)^8$.

 iii) $X_4$ (which is not isolated) is isomorphic to the Jacobian
 variety of the Klein's smooth hypersurface of degree 3 in $\Bbb
 P^4$. \endproclaim

 \demo{Proof} The first assertion follows from an explicit
 computation: only $C_1$, $C_2$ and $C_3$ have trivial isotropy
 group and we can apply Theorem 1.

 Assertions i) and ii) follows
 from the computation of explicit basis of $H^0(Y,\omega_Y)$
 ($Y$ denoting the corresponding curve). In
 these spaces $\Bbb Z_{11}$ acts in a natural way, it is a simple
 calculus to show that the analytic spectrum of these action
 correspond, respectively to $C_1$ and $C_2$. For instance, for i), a basis
 for $H^0(Y, \omega_Y)$ if $Y$ is given by the above equation is

 $$\{x^{i-1}\frac{dx}{y} \mid i=1,...,g\},$$
 a simple computation shows that $\Bbb Z_{11}$ acts on this basis by permutation, and the induced
 representation is given by a diagonal matrix with set of
 eigenvalues equal to $C_1$. Analogously, a basis for the space of
 differentials on the Lefschetz curve can be obtained from
 adjunction, and a similar calculus on the action of $\Bbb Z_{11}$ leads to the
 desired conclusion.

 Finally, part iii) is a consequence of [1]. Note that in this
 case $X_4$ admits $PSL(2,\Bbb F_{11})$ as the complete group of
 automorphism (compare with [2], Thm. 2).

 \enddemo

 We remark that $X_4$ gives another example of an inclusion of the
 type $\Cal A_g(p,\alpha)\subset \Cal A_g(q,\beta)$, this time with $g=5$,
 $p=11$ and $q=5$.

\enddocument